\documentclass{article}
\usepackage[usenames,dvipsnames]{pstricks}

\usepackage{comment}
\usepackage{pstricks,psfrag}
\usepackage{xcolor,rotating,epic,eepic}
\usepackage{amssymb}
\usepackage{amsmath}
\usepackage{times} 
\usepackage{epsfig,t1enc}
\usepackage{graphicx}
\usepackage{graphics}
\newtheorem{defi}{Definition}
\newtheorem{theorem}{Theorem}
\newtheorem{pro}{Proposition}

\def\II{\hbox{1\kern-.2em\hbox{I}}}

\date{}
\title{An iterative approach for real roots of polynomials} 

\author{J.-M. Billiot \and E. Fontenas
\thanks{Laboratoire Jean Kuntzmann, UMR 5224, Universit\'e Grenoble Alpes, 700 avenue centrale,
38041 Domaine Universitaire de Saint-Martin-d'H\`eres, France.
Jean-Michel.Billiot@univ-grenoble-alpes.fr, Eric.Fontenas@univ-grenoble-alpes.fr}}

\begin{document}
\maketitle
\begin{abstract}
In the present study, we propose  necessary and sufficient assumptions on the coefficients in order to only get distinct real roots of polynomials.  
\end{abstract}
\vspace{1cm}

\par 
{\em AMS}: 26A06, 26C10.
\section{Introduction}

It is well known that polynomials are very useful in order to approximate functions. The research of roots of polynomials is an old and famous problem. The theory of equations was studied by prestigious mathematicians such as d'Alembert, Cauchy, Gauss, Euler, Lagrange, Hermite, Galois among others. Polynomials with only real zeros arise often in different branches of mathematics. Nowadays, this is always the subject of an intense research, for example see \cite{Liu07} and references therein.

\cite{Liu07} proposed a unified approach to
polynomial sequences with only real zeros.
 They give new sufficient conditions for a sequence of polynomials to have only real zeros
based on the method of interlacing zeros. As applications, they derived the reality of zeros of orthogonal polynomials, matching polynomials, Narayana
polynomials and Eulerian polynomials.

Recently \cite{Gonz18}, \cite{Gonz19} studied cubic, quartic and quintic
polynomials and proposed conditions on the coefficients derived  from the Sturm sequence that will determine the real and complex root multiplicities together with the order of the real roots with respect to multiplicity.

Historically, equations of the first and second degree (where the coefficients are given numbers) are already solved with a general method by the Babylonians around 1700 bc. J.C and may be even earlier.

For the equations of degree three, it is necessary to wait until 1515 with the Italian Scipio del Ferro (1465-1526) whose papers are however lost.
Then, his compatriots Nicolo Tartaglia and G\'erolamo Cardano (1501-1576) continue his work. But it is Euler (1707-1783) who clarified the determination of the three roots in a Latin article of 1732.

For the equations of the degree four, the answer comes from Jerome Cardan (1501-1576) and Lodovico Ferrari (1522-1565). Cardan gives a method in chapter 39 of the Ars Magna. He states that it was found by his pupil Lodovico Ferrari.
    In 1615, Fran\c cois Vi\`ete (1540-1603) clearly explains Ferrari's method. Descartes (1596-1650) also exposes another method of resolution.

For the equations of  fifth order  and above, the theorem, sometimes called Abel-Ruffini's theorem, indicates that:
    "For every polynomial with coefficients of degree greater than or equal to five, there is no expression by radicals of the roots of the polynomial, that is to say of expression using only the coefficients, the value one, the four operations and the extraction of the nth roots".

This result is expressed for the first time by Paolo Ruffini (1765-1822), then rigorously proved by Niels Henrik Abel (1802-1829).

However, it is the legendary French mathematician Evariste Galois (1811-1832) who gives a necessary and sufficient condition for a polynomial equation to be solvable by radicals. He introduces permutations groups of the roots, now called Galois groups.
 This more precise version makes it possible to exhibit equations of degree five, with integer coefficients, whose complex roots which exist according with D'Alembert-Gauss's theorem do not express themselves by radicals.

In this paper, we propose an iterative approach for real roots of polynomials based on the ideas of C Sturm \cite{Stu35}. First of all, we give the family of polynomials we are interested in. For simplicity reason (but our method works also in the general case), we are interested in polynomials having the coefficient  of highest power equal one (if not we have just to divide the polynomial by this coefficient). Besides, we assume the coefficient of the second highest power is zero: we can always produce a translation that leads to this form. This means that the sum of the roots vanishes. We take as example the well known case of the order three. Then, we explain the general idea.\\ 
The goal is to give the motivation and convince the reader what kind of results we prove. Necessary and sufficient assumptions on the coefficients are given in order to obtain only distinct real roots. More precisely, we build a simple characterization of interlaced roots of the remainder of the polynomial and its derivative using  extrema of the starting polynomial. We distinguish two cases according with the degree even or odd of the studied polynomial.
 In fact, the originality of our method is that we just need the same assumption for all degrees of the starting polynomial. 
 
After recalling a consequence of the Sturm's theorem, we can identify the greatest common divisor (GCD) of a polynomial and its derivative with  resultants and Sylvesters' matrices.
 But even if, as describe in \cite{VanVl99} and \cite{Akr88}, studying remainders of the Sturm sequence can be seen as  minors of a single determinant, it remains in general a very difficult and laborious task.  
 
We  see in a new light the real roots of polynomials of orders three, four, five, six and seven.
We use maple software to calculate different remainders of Euclidean divisions. We interpret our result because, for low degree, the expression of the roots are available since a long time (for example for orders three and four.) We think also interesting to compare our result with the assumptions arising from the Sturm's approach.\\
Besides, we concentrated on particular cases of multiple roots. This is because in such cases, the intervals we are given are reduced to a point. As well, the same limit cases are obtained using the Sturm's approach. On the other hand, it shows a way to establish upper and lower bounds for these intervals.\\  
The case of the fifth order is studied in details and we explain how our assumptions can be easily satisfied.  In some sense, our result can be considered as an extension of the Sturm's theorem and  wonderful ideas he exhibited in \cite{Stu35}. We finish with some concluding remarks and perspectives. 
\section{Method}
In this section, we introduce the family of studied  polynomials. As example, we present the case of the order three.
Next, we give the general idea. In fact, we will precise this later when we will present the cases of order five, six and seven. 
The idea here is to build a series of polynomials all having real roots introducing iteratively the same assumption on the last constant of the built polynomial. We consider the following polynomial sequence:
$$\left\{\begin{array}{l}
 P_2(x)=x^2+c_0/3\\
 P_3(x)= x^3+c_0x+c_1\\
 P_4(x)= x^4+\frac{4!}{ 2!3!}c_0x^2+\frac{4!}{1!3!}c_1x+\frac{4!}{3!0!}c_2\\
P_5(x)=x^5+\frac{5!}{3!3!}c_0x^3+\frac{5!}{3!2!}c_1x^2+\frac{5!}{1!3!}c_2x+\frac{5!}{0!3!}c_3\\
\cdots\\
 P_n(x)=x^n+a_{n-2}x^{n-2}+a_{n-3}x^{n-3}+\cdots+a_0
\end{array}\right.$$
where
$$\forall k\in\{0,\ldots,n-2\},\quad a_k=\frac{n!}{k!3!}c_{n-2-k}.$$
We note that $$\forall n\geq 3,\, P_{n-1}(x)=\frac{1}{n}P^{'}_{n}(x).$$
As the notion of interlacing is crucial, a definition is welcome:
\begin{defi}
Given both polynomials $P$ and $Q$ of order $n$ and $n-1$ respectively and  $\{\alpha_i\}_{1\leq i\leq n}$ and $\{\beta_j\}_{1\leq j\leq n-1}$ be all real roots of $P$ and $Q$ in nonincreasing order respectively. We say that the roots of $Q$ are interlaced with the roots of $P$ if 
$$\alpha_1\leq \beta_1\leq \alpha_2\leq \beta_2\leq \cdots \leq \beta_{n-2}\leq \alpha_{n-1}\leq \beta_{n-1}\leq \alpha_{n}.$$
\end{defi}
\subsection{Examples of polynomials of degree two and three}
\noindent 1. Consider $P_2(x)=x^2+p/3$. This polynomial has two distinct real roots as soon as $p<0$.\\
\noindent 2. For the polynomial $P_3(x)=x^3+px+q$, by the Euclidean division by $P'_3=3P_2$, it comes then 
$$P_{3}(x)=\frac{x}{3}P^{'}_{3}(x)-R_{1}(x)=xP_{2}(x)- R_{1}(x)$$
with $R_{1}(x)=-(\frac{2}{3}px+q)$. The polynomial $P_3$  has three real roots if and only if $P_2$  has two real roots and if $R_1$ has a real root interlaced with those of  $P_2$. Because of the previous remark, $P_2$  has two real roots if $p<0$.\\
It then remains to fix the constant $q$ so that the root of $R_1$ is interlaced with those of $P_2$. The root of $R_{1}$ is $\beta^{(1)}_{1}=-\frac{3q}{2p}$ and is interlaced with those of $P_{2}$ if $$-\frac{3q}{2p}\in \left]\alpha^{(2)}_{1}=-\left(-\frac{p}{3}\right)^{1/2};\alpha^{(2)}_{2}=\left(-\frac{p}{3}\right)^{1/2}\right[$$ or $$q\in \left]-\frac{2}{3}p\alpha^{(2)}_{1}; -\frac{2}{3}p\alpha^{(2)}_{2}\right[.$$\\
If we call $R_{1}^{0}$ the function defined by  $R_{1}^{0}(x)=R_1(x)+q=-\frac{2}{3}px$,  we have that $$q\in]R_{1}^{0}(\alpha^{(2)}_{1}); R_{1}^{0}(\alpha^{(2)}_{2})[.$$
This means that $R_{1}\left(\alpha^{(2)}_{1}\right)<0$ and that $R_{1}\left(\alpha^{(2)}_{2}\right)>0$. Note here that the result means that the discriminant of $P_3$ namely $\Delta(P_{3})= -(4p^3+27q^2)$ is not negative.
Moreover, if $p<0$, it is well known that only cases of multiple root for $P_{3}$ correspond to $\beta^{(1)}_{1}= \alpha^{(2)}_{1}$ or $\beta^{(1)}_{1}=\alpha^{(2)}_{2}$ according with the sign of $q$. Remark that this real double root is also a root of $P_{2}$ and $R_{1}$.

\subsection{General idea}
Consider the  family of polynomials defined previously

$$P_n(x)=x^n+a_{n-2}x^{n-2}+a_{n-3}x^{n-3}+\cdots+a_0$$
where
$$\forall k\in[0;n-2],\quad a_k=\frac{n!}{k!3!}c_{n-2-k}.$$

\noindent From a more general point of view, the result is based on the following remarks:\\

\noindent 1. Set $P_n(x)=xQ_{n-1}(x)-R_{n-2}(x)$.
If $Q_{n-1}$ has $n-1$ real roots and if $R_{n-2}$ has $n-2$ real roots interlaced with those of $Q_{n-1}$, then $P_n$ has $n$ real roots.\\
2. For a polynomial of degree $n$, to have $n$ real roots, it is necessary that its derivative has
$n-1$ distinct real roots $\alpha^{(n-1)}_{1}\ldots\alpha^{(n-1)}_{n-1}$. So,
$P_n(x)=\frac{x}{n}P'_{n}(x)-R_{n-2}(x)$ has $n$ distinct real roots as soon as its
derivative has $n-1$ distinct real roots and those of $R_{n-2}$ denoted by $\beta_{1}^{(n-2)}\ldots \beta_{n-2}^{(n-2)}$ are interlaced  with those of $P'_{n}$.\\

The coefficient of highest degree of $R_{n-2}$ is $-\frac{2}{n}a_{n-2}$. If $a_{n-2}<0$ (as the sign of the discriminant of $P_{2}$) and if $R_{n-2}$ has $n-2$  real distinct and interlaced roots with $P^{'}_{n}$ then $P_{n}$ has $n$ distinct real roots.
If we write:
$$R_{n-2}(x)=\frac{1}{n-2}(x-\beta^{(1)}_{1})R^{'}_{n-2}(x)-T_{n-4}(x),$$ as $R^{'}_{n-2}(x)=(n-1)R_{n-3}(x)$, $R^{'}_{n-2}$ has $n-3$ distinct real roots. It is therefore necessary to ensure that $T_{n-4}$ has $n-4$ distinct real roots denoted by $\gamma^{(n-4)}_{1}\ldots \gamma^{(n-4)}_{n-4}$ and besides interlaced with those  $\beta_{1}^{(n-3)}\ldots \beta_{n-3}^{(n-3)}$ of $R^{'}_{n-2}$. The highest degree coefficient
of $T_{n-4}$  is at a positive factor proportional to $\Delta_{2}$ (discriminant of $R_{2}$). Then
$$T_{n-4}(x)=\frac{1}{n-4}\left(x-\gamma^{(1)}_{1}\right)T^{'}_{n-4}(x)-U_{n-6}(x).$$ 
The roots of $U_{n-6}$ noted $\delta^{(n-6)}_{1}\ldots \delta^{(n-6)}_{n-6}$ will be interlaced with those $\gamma_{1}^{(n-5)}\ldots \gamma_{n-5}^{(n-5)}$ of $T^{'}_{n-4}$. The coefficient of highest degree of $U_{n-6}$ is a positive factor proportional to the
discriminant of $T_{2}$.  Then we have:
$$U_{n-6}(x)=\frac{1}{n-6}\left(x-\delta^{(1)}_{1}\right)U^{'}_{n-6}(x)-V_{n-8}(x)$$ 
and so on.
We will specify this construction later when studying polynomials of degree five, six and seven, when we will explain how  choosing iteratively the last coefficient $a_{0}$ of the polynomial $P_{n}$. 
This leads us now to the presentation of our main result.

\section{Main result}
First of all, we present a characteristic property of interlaced roots. Then, we explain how choosing the integration constant $ a_ {0} $ using the extrema of  $P_{n}$ and the remainder of Euclidean division of $P_{n}$ by $P^{'}_{n}$. At last, we propose some necessary and sufficient assumptions  on the coefficients in order to obtain $n$ distinct real roots for $P_{n}$.
\begin{pro}
Let $P_{n}$ be a polynomial and let $$P_{n}(x)=\frac{x}{n}P'_{n}(x)-R_{n-2}(x)=xP_{n-1}(x)-R_{n-2}(x).$$
Denote by $\alpha^{(n-1)}_{1}\ldots\alpha^{(n-1)}_{n-1}$, the $n-1$ distinct real roots of  $P_{n-1}$ and  $\beta^{(n-2)}_{1}\ldots\beta^{(n-2)}_{m}$, the $n-2$ distinct real roots of  $R_{n-2}$.\\
If $n$ is even, $$\alpha_{1}^{(n-1)}<\beta_{1}^{(n-2)}<\alpha_{2}^{(n-1)}<\beta_{2}^{(n-2)}<\ldots <\alpha_{n-2}^{(n-1)}<\beta_{n-2}^{(n-2)}<\alpha_{n-1}^{(n-1)}$$
$$ \Leftrightarrow \sup_{k\in{\{1,\ldots, \frac{n-2}{2}}\}}R_{n-2}(\alpha_{2k}^{(n-1)})<0<\inf_{k\in{\{0,\ldots, \frac{n-2}{2}}\}}R_{n-2}(\alpha_{2k+1}^{(n-1)}).$$\\
If $n$ is odd,
$$\alpha_{1}^{(n-1)}<\beta_{1}^{(n-2)}<\alpha_{2}^{(n-1)}<\beta_{2}^{(n-2)}<\ldots< \alpha_{n-2}^{(n-1)}<\beta_{n-2}^{(n-2)}<\alpha_{n-1}^{(n-1)}$$
$$ \Leftrightarrow \sup_{k\in{\{0,\ldots, \frac{n-3}{2}}\}}R_{n-2}(\alpha_{2k+1}^{(n-1)})<0<\inf_{k\in{\{1,\ldots, \frac{n-1}{2}}\}}R_{n-2}(\alpha_{2k}^{(n-1)}).$$

\end{pro}
\noindent{\bf Proof:} We only present here the proof for $n $
even.
Suppose that $R_{n-2}$ admits $n-2$ distinct roots $\beta_{1}^{(n-2)}\ldots \beta_{n-2}^{(n-2)}$ interlaced with those of $P_{n-1}$ : 
$$\alpha_{1}^{(n-1)}<\beta_{1}^{(n-2)}<\alpha_{2}^{(n-1)}<\beta_{2}^{(n-2)}<\ldots <\alpha_{n-2}^{(n-1)}<\beta_{n-2}^{(n-2)}<\alpha_{n-1}^{(n-1)}.$$
So, 
$$\sup_{k\in{\{0,\ldots, \frac{n-2}{2}}\}}P_{n}\left(\alpha_{2k+1}^{(n-1)}\right)<0<\inf_{k\in{\{1,\ldots, \frac{n-2}{2}}\}}P_{n}\left(\alpha_{2k}^{(n-1)}\right).$$
which is also written
$$\sup_{k\in{\{1,\ldots, \frac{n-2}{2}}\}}R_{n-2}\left(\alpha_{2k}^{(n-1)}\right)<0<\inf_{k\in{\{0,\ldots, \frac{n-2}{2}}\}}R_{n-2}\left(\alpha_{2k+1}^{(n-1)}\right).$$
Conversely, 

$$\begin{array}{c}
\displaystyle\sup_{k\in{\{1,\ldots, \frac{n-2}{2}}\}}R_{n-2}\left(\alpha_{2k}^{(n-1)}\right)<0<\inf_{k\in{\{0,\ldots, \frac{n-2}{2}}\}}R_{n-2}\left(\alpha_{2k+1}^{(n-1)}\right)\\
\Longleftrightarrow \\
\forall k\in\{1,..n-2\},\,R_{n-2}\left(\alpha_{k}^{(n-1)}\right)\times R_{n-2}\left(\alpha_{k+1}^{(n-1)}\right)<0
\end{array}$$
This implies that $R_{n-2}$ has $n-2$ interlaced roots with those of $P_{n-1}$. \\

The following theorem allows to choose $a_0$ so that the new polynomial $P_n$ has $n$ distinct real roots:
\begin{theorem}
Let $P_n$ be a polynomial and set $$P_{n}(x)=\frac{x}{n}P^{'}_{n}(x)-R_{n-2}(x)=xP_{n-1}(x)-R_{n-2}(x).$$
We call $R_{n-2}^0=R_{n-2}+a_0$.\\
Denote by $\alpha^{(n-1)}_{1}\ldots\alpha^{(n-1)}_{n-1}$  the $n-1$ distinct real roots of  $P_{n-1}$.\\

\noindent For $n$ even, if
$$\left\{\begin{array}{l}
\displaystyle\sup_{k\in{\{1,\ldots, \frac{n-2}{2}}\}}R_{n-2}(\alpha_{2k}^{(n-1)})<0<\displaystyle\inf_{k\in{\{0,\ldots, \frac{n-2}{2}}\}}R_{n-2}(\alpha_{2k+1}^{(n-1)})\\
a_{0}\in\left]\displaystyle\sup_{k\in{\{1,\ldots, \frac{n-2}{2}}\}}R_{n-2}^{0}(\alpha_{2k}^{(n-1)});\displaystyle\inf_{k\in{\{0,\ldots, \frac{n-2}{2}}\}}R_{n-2}^{0}(\alpha_{2k+1}^{(n-1)})\right[,
\end{array}\right.$$
then $P_{n}$ has $n$ distinct real roots.\\

\noindent For $n$ odd, if
$$\left\{\begin{array}{l}
\displaystyle\sup_{k\in{\{0,\ldots, \frac{n-3}{2}}\}}R_{n-2}(\alpha_{2k+1}^{(n-1)})<0<\displaystyle\inf_{k\in{\{1,\ldots, \frac{n-1}{2}}\}}R_{n-2}(\alpha_{2k}^{(n-1)})\\
a_0\in\left]\displaystyle\sup_{k\in{\{0,\ldots, \frac{n-3}{2}}\}}R^{0}_{n-2}(\alpha_{2k+1}^{(n-1)})\,;\,\displaystyle\inf_{k\in{\{1,\ldots, \frac{n-1}{2}}\}}R^{0}_{n-2}(\alpha_{2k}^{(n-1)})\right[,
\end{array}\right.$$
then $P_{n}$ has $n$ distinct real roots.
\end{theorem}
\noindent{\bf Remarks:} 1. For example, in the case of $n$ odd, if $a_0=\displaystyle\sup_{k\in{\{0,\ldots, \frac{n-3}{2}}\}}R^{0}_{n-2}(\alpha_{2k+1}^{(n-1)})$ or $a_0=\displaystyle\inf_{k\in{\{1,\ldots, \frac{n-1}{2}}\}}R^{0}_{n-2}(\alpha_{2k}^{(n-1)})$, that means that the polynomial $P_n$ as a double root which is a root of $P_{n-1}$ and $R_{n-2}$.\\
2. Either $a_{0}=\displaystyle\sup_{k\in{\{0,\ldots, \frac{n-3}{2}}\}}R^{0}_{n-2}(\alpha_{2k+1}^{(n-1)})=\displaystyle\inf_{k\in{\{1,\ldots, \frac{n-1}{2}}\}}R^{0}_{n-2}(\alpha_{2k}^{(n-1)})$ reach in two real distinct roots of $P_{n-1}$: $\alpha^{(n-1)}_{i}$ and $\alpha^{(n-1)}_{j}$.
So these two real distinct roots are also
roots of $R_{n-2}$ and therefore of $P_n$ .
We deduce that these roots are double roots of $P_n$.\\
Or, either, $a_{0}=\displaystyle\sup_{k\in{\{0,\ldots, \frac{n-3}{2}}\}}R^{0}_{n-2}(\alpha_{2k+1}^{(n-1)})=\displaystyle\inf_{k\in{\{1,\ldots, \frac{n-1}{2}}\}}R^{0}_{n-2}(\alpha_{2k}^{(n-1)})$ reached in a real double roots of $P_{n-1}$ then this is a triple root of $P_{n}$ and a double root of $R_{n-2}$. \\

Applying recursively Theorem 1, the following theorem allows to choose all coefficients $a_l$, $0\leq l<n-2,$ 
\begin{theorem}
Let $P_{n}(x)=x^{n}+a_{n-2}x^{n-2}+\ldots+ a_{1}x+a_{0}$ and define the sequence $[P_{i},R_{i},R^{0}_{i}]$, $i\in\{3,\ldots,n\}$, such that 
$$\left\{\begin{array}{l}
P_{i-1}(x)=\displaystyle\frac{1}{i}P'_{i}(x)\\
P_{i}(x)=xP_{i-1}(x)-R_{i-2}(x)\\
\end{array}\right.$$
and, $\forall i\in \{3,\ldots,n\},\, R^{0}_{i-2}(x)=R_{i-2}(x)+a_{n-i}$.\\
If $n$ is even, $P_n$ has  $n$ distinct real roots if and only if, for all $l\in\{0,...,\frac{n}2-2\}$,
$$\left\{\begin{array}{l}
\bullet \, a_{n-2}<0\\
\bullet \, \displaystyle\sup_{k\in \{1,\ldots,\frac{n-2}{2}-l\}}R^{0}_{n-2l-2}(\alpha^{(n-2l-1)}_{2k})\,<\, \inf_{k\in \{0\ldots \frac{n-2}{2}-1\}}R^{0}_{n-2l-2}(\alpha^{(n-2l-1)}_{2k+1})\\
\bullet\, a_{2l}\in \left]\displaystyle\sup_{k\in \{1,\ldots, \frac{n-2}{2}-l\}}R^{0}_{n-2l-2}(\alpha^{(n-2l-1)}_{2k})\,;\, \inf_{k\in \{0,\ldots ,\frac{n-2}{2}-1\}}R^{0}_{n-2l-2}(\alpha^{(n-2l-1)}_{2k+1})\right[\\
\bullet \,\displaystyle\sup_{k\in \{0,\ldots, \frac{n-4}{2}-l\}}R^{0}_{n-2l-3}(\alpha^{(n-2l-2)}_{2k+1})\,<\,  \inf_{k\in \{1,\ldots, \frac{n-2}{2}-l\}}R^{0}_{n-2l-3}(\alpha^{(n-2l-2)}_{2k})\\
\bullet \,a_{2l+1}\in\left] \displaystyle\sup_{k\in \{0,\ldots, \frac{n-4}{2}-l\}}R^{0}_{n-2l-3}(\alpha^{(n-2l-2)}_{2k+1})\,;\,  \inf_{k\in \{1,\ldots, \frac{n-2}{2}-1\}}R^{0}_{n-2l-3}(\alpha^{(n-2l-2)}_{2k})\right[.
\end{array}\right.$$
If $n$ is odd, $P_n$ has  $n$  distinct real roots if and only if,\\
$\bullet$ $a_{n-2}<0$.\\
$\bullet$ $\forall\, l\in\{0,...,\frac{n-3}2\},\,
\left\{\begin{array}{l}\bullet \displaystyle\sup_{k\in \{0,\ldots, \frac{n-3}{2}-l\}}R^{0}_{n-2l-2}(\alpha^{(n-2l-1)}_{2k+1})\,<\, \inf_{k\in \{0,\ldots, \frac{n-1}{2}-l\}}R^{0}_{n-2l-2}(\alpha^{(n-2l-1)}_{2k})\\
\bullet\, a_{2l}\in\left] \displaystyle\sup_{k\in \{0,\ldots, \frac{n-3}{2}-l\}}R^{0}_{n-2l-2}(\alpha^{(n-2l-1)}_{2k+1})\,;\, \inf_{k\in \{0,\ldots, \frac{n-1}{2}-l\}}R^{0}_{n-2l-2}(\alpha^{(n-2l-1)}_{2k})\right[
\end{array}\right.$\\
$\bullet$ $\forall\, l\in\{0,...,\frac{n-5}2\},\,\left\{\begin{array}{l}\bullet \displaystyle\sup_{k\in \{1,\ldots, \frac{n-3}{2}-l\}}R^{0}_{n-2l-3}(\alpha^{(n-2l-2)}_{2k})\,<\, \inf_{k\in \{0,\ldots, \frac{n-3}{2}-1\}}R^{0}_{n-2l-3}(\alpha^{(n-2l-2)}_{2k+1})\\
\bullet\, a_{2l+1}\in\left]\displaystyle\sup_{k\in \{1,\ldots, \frac{n-3}{2}-l\}}R^{0}_{n-2l-3}(\alpha^{(n-2l-2)}_{2k})\,;\, \inf_{k\in \{0,\ldots, \frac{n-3}{2}-1\}}R^{0}_{n-2l-3}(\alpha^{(n-2l-2)}_{2k+1})\right[.
\end{array}\right.$
\end{theorem}

\section{Comparison with Sturm's approach}
Of course, at this stage, it is difficult to see what our result means. An interpretation of our assumptions should be welcome. In particular, we may ask if we can compare our assumptions to
those arising from Sturm's theorem.
We denote $S_{n}$ the first term of the Sturm sequence of the polynomial $P_{n}$ defined as follows:

$$\left\{\begin{array}{l}
S_{n}=P_{n}\\
S_{n-1}=P^{'}_{n}\\
S_{n}=Q_{1}S_{n-1}-S_{n-2}\\
\ldots\\
S_{2}=Q_{n-1}S_{1}-S_{0}.
\end{array}\right.$$

\noindent And, with our notations, 

$$\left\{\begin{array}{l}
P_{n}(x)=xP_{n-1}(x)-R_{n-2}(x)\\
R_{n-2}(x)=R^{0}_{n-2}(x)-a_{0}.
\end{array}\right.$$
\begin{pro}
$$
\begin{array}{lll}
S_{0}&=&K_{1}\Delta(P_{n})=K_{2}\displaystyle\prod_{i=1}^{n-1}P_{n}(\alpha^{(n-1)}_{i})
\\
&=&(-1)^{n-1}K_{2}\displaystyle\prod_{i=1}^{n-1}R_{n-2}(\alpha^{(n-1)}_{i})=(-1)^{n-1}K_{2}\prod_{i=1}^{n-1}[R^{0}_{n-2}(\alpha_{i}^{(n-1)})-a_0]\end{array}$$
 where $\alpha^{(n-1)}_{i}, i=1,\ldots, n-1$, distinct real roots of $P^{'}_{n}$ and $K_{1}$, $K_{2}$ strictly not negative constants. $\Delta(P_{n})$ represents the discriminant of $P_{n}$.
\end{pro}
\noindent{\bf Remarks:}\\
1. It is well known that, if the discriminant of the polynomial $P_n$ is not negative, it is a necessary but not sufficient assumption to obtain only real roots. The previous proposition follows from well known results on resultants and Sylvesters' matrices see for example \cite{Gelf94}.\\
2. A direct consequence of Sturm's theorem is that $P_{n}$ has $n$ real roots if the terms of
higher degree of $S_{j},\,j\in\{0,\ldots,n\}$, are all not negative.\\

For illustrative purposes, let us describe what happen for polynomials of degree three, four, five, six and seven. Everytime, we study different cases of multiple roots up to order six.

\subsection{Polynomial of order three}
In our case, the polynomial $P_3$ has three roots if $p<0$ and $4p^3+27q^2<0$. By Sturm's method,
$$\left\{\begin{array}{l}
S_{3}(x)=P_{3}(x)=x^3+px+q\\
S_{2}(x)=P^{'}_{3}(x)=3x^2+p\\
S_{1}(x)=R_{1}(x)=-(\frac{2}{3}px+q)\\
S_{0}=K_{2}\displaystyle\prod_{i=1}^{2}P_{3}(\alpha_{i}^{(2)})= K_{2}\displaystyle\prod_{i=1}^{2}R_{1
}(\alpha_{i}^{(2)})=
K_{2}\displaystyle\prod_{i=1}^{2}(R^{0}_{1}(\alpha_{i}^{(2)})-q)=K_{2}\left(-\frac{4p^3}{27}-q^2\right)=K_{1}\Delta(P_{3})
\end{array}\right.$$
The form of $S_{0}$ can be deduced easily from proposition $2$.

\subsection{Polynomial of order four}
The Sturm's sequence is given by
$$\left\{\begin{array}{l}
S_{4}(x)=P_{4}(x)=x^4+2px^2+4qx+4r\\
S_{3}(x)=P^{'}_{4}(x)=4x^3+4px+4q\\
S_{2}(x)=R_{2}(x)=-px^2-3qx-4r\\
S_{1}(x)=\displaystyle\frac{-1}{p^2}[(-4pr+p^3+9q^2)x +q(12r+p^2)]\\

\end{array}\right.$$
and 

$$\begin{array}{lll}S_{0}&=&\displaystyle\frac{p^2}{(-4pr+p^3+9q^2)^2}(64r^3-32p^2r^2+4p^4r+72prq^2-27q^4-2p^3q^2)\\
&=&\displaystyle\frac{p^2\Delta(P_{4})}{256(-4pr+p^3+9q^2)^2}\\
&=& K_{2}\displaystyle\prod_{i=1}^{3}P_{4}(\alpha_{i}^{(3)})=-K_{2}\prod_{i=1}^{3}R_{2}(\alpha_{i}^{(3)})=-K_{2}\prod_{i=1}^{3}(R^{0}_{2}(\alpha_{i}^{(3)})-4r)
\end{array}$$
with  $\alpha_{i}^{(3)},i=1..3$, the three distinct roots of $P_3$,  $\Delta(P_4)$ the discriminant of $P_4$ and $R_{2}(x)= -px^2-3qx-4r$,  $R^{0}_{2}(x)= -px^2-3qx$.\\
Using Sturm's theorem, we have four distinct real roots if the coefficients of the term of highest
degree of $S_{3}, S_2, S_1$ and $S_0$ are strictly positive. These assumptions are the following: $$p<0,\,-4pr+p^3+9q^2<0,\,\Delta(P_{4})>0.$$
Using our method, we have four distinct real roots for $P_{4}$ if we choose:
$$p<0,\quad q\in \left]-2(-\frac{p}{3})^{3/2}; 2(-\frac{p}{3})^{3/2}\right[,\quad 4r\in ]R_{2}^{0}(\alpha_{2}^{(3)});\inf_{i\in{\{1,3\}}}R^{0}_{2}(\alpha_{i}^{(3)}) [\subset ]R^{0}_{2}(\beta^{(1)}_{1});+\infty[ $$
where $\beta^{(1)}_{1}=-\frac{3q}{2p}$ is the root of $R_{1}$.\\
\begin{enumerate}
\item Say $4r\in ]R^{0}_{2}(\beta^{(1)}_{1});+\infty[$ means $R_2$ has two distinct roots $\beta_{1}^{(2)}$ and $\beta_{2}^{(2)}$ (the discriminant of $R_2$ is given by $\Delta_{2}=9q^2-16pr>0$). Clearly:
$$\Delta_{2}>0\Leftrightarrow R_{2}(\beta^{(1)}_{1})=R^{0}_{2}(\beta^{(1)}_{1})-4r<0.$$
\item Say $4r\in ]R_{2}^{0}(\alpha_{2}^{(3)});\displaystyle\inf_{i\in{\{1,3\}}}R^{0}_{2}(\alpha_{i}^{(3)}) [$ means the roots of $R_2$ are interlaced with $\alpha_{1}^{(3)}$, $\alpha_{2}^{(3)}$, $\alpha_{3}^{(3)}$.
\item If $q<0$, we get:
$\beta_{1}^{(2)}$ and $\beta_{2}^{(2)}$ interlaced with $\alpha_{1}^{(3)}, \alpha_{2}^{(3)},\alpha_{3}^{(3)}$ is equivalent to $4r\in ]R_{2}^{0}(\alpha_{2}^{(3)}); R_{2}^{0}(\alpha_{1}^{(3)})[$.  Moreover we can show as $\alpha_{3}^{(3)}> 3q/p$ that $\alpha_{2}^{(3)}$ is 
closer to $-\frac{3q}{2p}$ than $\alpha_{1}^{(3)}$ and so that $R_{2}^{0}(\alpha_{2}^{(3)})< R_{2}^{0}(\alpha_{1}^{(3)})$. This means $R_{2}(\alpha_{2}^{(3)})<0$ and 
$R_{2}(\alpha_{1}^{(3)})>0$. In that case, we remark that $\alpha^{(3)}_{1}<0$, $\alpha^{(3)}_{2}<0$ and $\alpha^{(3)}_{3}>0$. So, $4r\in ]R_{2}^{0}(\alpha_{2}^{(3)});\displaystyle\inf_{i\in{\{1,3\}}}R^{0}_{2}(\alpha_{i}^{(3)}) [=]R_{2}^{0}(\alpha_{2}^{(3)});R^{0}_{2}(\alpha_{1}^{(3)}) [$.
\item If $q>0$, we have
$4r\in ]R_{2}^{0}(\alpha_{2}^{(3)});\displaystyle\inf_{i\in{\{1,3\}}}R^{0}_{2}(\alpha_{i}^{(3)}) [=]R_{2}^{0}(\alpha_{2}^{(3)});R^{0}_{2}(\alpha_{3}^{(3)}) [$.

\end{enumerate}

\noindent{\bf Study of multiple roots for a polynomial of order four}\\
1. If $r=-p^2/12$, we are in a limit case in the following sense 
$$\Delta(P_4)=-\frac{256}{27}(4p^3+27q^2)^2>0$$
$$-4pr+p^3+9q^2=\frac{1}{3}(4p^3+27q^2)<0.$$
This implies $4p^3+27q^2=0$. We deduce that $S_{1}$ is
identically zero.
According with the sign of $q$, $P_4$ has a triple root $\pm\sqrt{-p/3}$  which is a double root of $P_{3}$ and of $R_2$.\\
 
\noindent 2. If $q=0$, the bounds become :
$\Delta(P_4)=64r(r-p^2/4)^2$ and the other is$-4p(r-p^2/4)$. Taking $r=p^2/4$, $S_{1}$ is
identically zero. This case corresponds with two double roots $\pm\sqrt{-p}$ for $P_{4}$. We can establish that the quotient of the Euclidean division of $P_{4}$ by $(x-a)^2$ is $3a^2+2xa+2p+x^2$. This one has two real roots if its discriminant $-8p-8a^2$ is not negative. This requires that $$a\in ]-\sqrt{-p};\sqrt{-p}[$$ which are precisely the double roots previously obtained for $q=0$.\\
 
\noindent Both previous points make it possible to find bounds for $r$:
$$4r\in ]R_{2}^{0}(\alpha_{2}^{(3)});\displaystyle \inf_{k\in \{1,3\}} R_{2}^{0}(\alpha_{k}^{(3)})[\subset \left]-\frac{p^2}{3};p^2\right[.$$

\subsection{Polynomial of order five}
In this section, we concentrate on the order five.
After giving the Sturm's assumptions, we express three particular cases. Then, we describe the assumptions of our theorem. In particular, the choice of $s$ is discussed using the general idea. We explain how our assumption can be satisfied.
 Next, different cases of multiple roots are specified. We can remark once more that our assumptions and Sturms' assumptions lead to the same result.
\subsubsection{The Sturm's assumptions}
The Sturm's sequence is given by
$$\left\{\begin{array}{l}
S_5(x)=P_{5}(x)=x^5+\frac{10p}{3}x^3+10qx^2+20rx+20s\\
S_{4}(x)=P^{'}_{5}(x)=5x^4+10px^2+20qx+20r\\
S_{3}(x)=R_{3}(x)=-(\frac{4}{3}px^3+6qx^2+16rx+20s)\\
S_{2}(x)=-5\left[\frac{(8p^3+81q^2-48pr)x^2}{4p^2}+\frac{(-15ps+4p^2q+54qr)x}{p^2}+4r+\frac{135qs}{2p^2}\right]\\
S_1(x)=-a_{S_1}\,x -b_{S_1}
\end{array}\right.$$
with \\
$\left\{\begin{array}{lll}
a_{S_1}&=&-80p^4r-2106q^2pr+1056p^2r^2-3456r^3+240p^2qs+3240qsr\\
&&+40p^3q^2+729q^4-450ps^2\\
b_{S_1}&=&-120p^4s-1755spq^2+1560p^2rs-4320r^2s+40p^3qr+729q^3r\\
&&-864qpr^2+2025s^2q.
\end{array}\right.$\\

\noindent$\begin{array}{lll}
S_{0}&=&1800s^2p^5-3600p^4qsr+1600r^3p^4-27000p^3s^2r-600q^2r^2p^3+1200p^3q^3s\\&&+37125p^2s^2q^2-23040r^4p^2
+50400p^2r^2qs
-85050prq^3s+38880r^3q^2p\\
&&+108000pr^2s^2-101250ps^3q+182250rq^2s^2-10935r^2q^4+21870q^5s\\
&&-259200r^3qs+82944r^5+50625s^4\\
&=&K_{1}\Delta(P_{5})=K_{2}\displaystyle\prod_{i=1}^{4}P_{5}(\alpha^{(4)}_{i})
=K_{2}\displaystyle\prod_{i=1}^{4}R_{3}(\alpha^{(4)}_{i})=K_{2}\displaystyle\prod_{i=1}^{4}(R^{0}_{3}(\alpha^{(4)}_{i})-20s).
\end{array}$\\

The result of Sturm gives five distinct real roots for $P_{5}$ under the following assumptions:

\begin{equation}
p<0
\end{equation}
\begin{equation}\label{eq1}
8p^3+81q^2-48pr<0
\end{equation}
\begin{equation}\label{eq2}
a_{S_1}<0
\end{equation}
\begin{equation}
\prod_{i=1}^{4}(R^{0}_{3}(\alpha^{(4)}_{i})-20s)>0.
\end{equation}
The polynomial $a_{S_1}$ in $s$ in the inequality (\ref{eq2}) having  $-450p>0$ as the coefficient before $s^2$ has to be not positive: we must have that its discriminant
$$144p^4q^2-2484q^2p^2r+11664q^2r^2-160p^5r+2112p^3r^2-6912pr^3+1458pq^4>0.$$
Its roots in $r$ are:
$$r_{1}=\frac{p^2}{6}+\frac{27q^2}{16p},\, r_2=\frac{5p^2-\sqrt{25p^4-648q^2p}}{72},\,r_{3}= \frac{5p^2+\sqrt{25p^4-648q^2p}}{72}.$$
In fact, three particular cases are interesting to explain.\\

\noindent 1. The particular case $q=0$ deserves to be detailed. It comes $r_1=p^2/6$, $r_2=5p^2/36$ and $r_3=0$. If we take $r=\frac{5p^2}{36}$, it gives two double roots $-(-\frac{5p}{3})^{1/2}$ and $(-\frac{5p}{3})^{1/2}$ and $s=0$.\\

\noindent 2. If $r_{1}=r_{2}$ or $r_1=r_3$, then $4p^3+27q^2=0$ or $8p^3+27^2q^2=0$. When $4p^3+27q^2=0$, the particular case $r=r_1=\frac{-p^2}{12}$ is interesting: the polynomial $P_4$ has one triple root. The first three Sturm's assumptions become
$$\left\{\begin{array}{l}
p<0\\
8p^3+81q^2-48pr=3(4p^3+27q^2)<0\\
a_{S_1}=(4p^3+27q^2)^2-450p(s+\frac{pq}{30})^2<0.
\end{array}\right.$$
$s$ is necessary equal to $\frac{-pq}{30}$: according with the sign of $q$, $\pm\sqrt{-p/3}$ is the quadruple root of $P_5$.\\
If $8p^3+27^2q^2=0$, the particular case $r=r_1=\frac{4p^2}{27}$ is interesting: the first three Sturms' assumptions become
$$\left\{\begin{array}{l}
p<0\\
8p^3+81q^2-48pr=\frac{1}{9}(8p^3+729q^2)=0\\
a_{S_1}= -450p(s-\frac{4pq}{5})^2+\frac{1}{729}(8p^3+729q^2)^2.
\end{array}\right.$$
This situation shows a limit case where, if $q<0$, taking $q=-\frac{\sqrt{-8p^3}}{27}$, $s$ is necessary equal to $\frac{4pq}{5}$. Then, $P_{5}$ have a real triple root $\alpha^{(5)}_{1}=\alpha^{(4)}_{1}=\alpha^{(3)}_{1}=-\frac{2\sqrt{-2p}}{3}$  and a real double root $\alpha^{(5)}_2=\alpha^{(4)}_2=\sqrt{-2p}$.
Notice that we have $$4r=-R_2^0(\alpha^{(3)}_{1})=\frac{16 p^2}{27}.$$ 
Using Sturms' polynomials, if $r=\frac{4p^2}{27}$ and $729q^2+8p^3=0$ and $s=\frac{4pq}{5}$ or, if $r=\frac{-p^2}{12}$ and $27q^2+4p^3=0$ and $s=\frac{2p^4}{405q}=\frac{-pq}{30}$, then $S_{2}$ is identically zero. Both
cases correspond respectively with a real triple root with a double real root and a  quadruple
real root and a simple real root for the polynomial $P_5$: this result is consistent with \cite{Gonz19}. 
\subsubsection{Our assumptions}
For our method, recall the assumptions of order four:
$$p<0,\quad q\in \left]-2\left(-\frac{p}{3}\right)^{3/2};2\left(-\frac{p}{3}\right)^{3/2}\right[,\quad 4r\in \left]R_{2}^{0}(\alpha_{2}^{(3)});\inf_{k\in{\{1,3\}}}R^{0}_{2}(\alpha_{k}^{(3)})\right[. $$
Applying the theorem 1, we explain how choosing the parameter $s$.\\
 
\noindent{\bf Choice of s}\\
First, $R_{3}$ must have three distinct real roots and
then they have to be interlaced with $\alpha^{(4)}_{{i}}, i=1\ldots 4$. According with the result of order three, we find
$$20s\in ]R^{0}_{3}(\beta^{(2)}_{2}),R^{0}_{3}(\beta^{(2)}_{1})[.$$
The Euclidean division of $R_{3}$  by $R^{'}_{3}$ 
give as a remainder whose sign we change
$$T_{1}(x)=\frac{2(16pr-9q^2)}{3p}x+20s-\frac{8qr}{p}$$
and this remainder vanishes at 
$$\gamma^{(1)}_{1}=\frac{3(-10sp+4qr)}{16pr-9q^2}$$
It is enough now that $\gamma^{(1)}_{1}$ is interlaced with
$\beta^{(2)}_{1}$ and $\beta^{(2)}_{2}$:
 $$\beta^{(2)}_{1}< \gamma^{(1)}_{1}< \beta^{(2)}_{2}.$$
We obtain :
$$ 20s\in ]R^{0}_{3}(\beta^{(2)}_{2});R^{0}_{3}(\beta^{(2)}_{1})[=\left] \frac{-9q^3}{p^2} +\frac{24qr}{p} -\frac{ (\Delta_{2})^{3/2}}{3p^2};\frac{-9q^3}{p^2} +\frac{24qr}{p} +\frac{ (\Delta_{2})^{3/2}}{3p^2}\right[.$$ 
With the assumption $\Delta_{2}=9q^2-16pr>0$ (the discriminant of $R_2$),  we deduce that  $R_{3}$ will have three distinct real roots.\\
Or, as $$-T_{1}(\beta^{(2)}_{1})=R_{3}(\beta^{(2)}_{1})= R^{0}_{3}(\beta^{(2)}_{1})-20s$$ and $$-T_{1}(\beta^{(2)}_{2})=R_{3}(\beta^{(2)}_{2})= R^{0}_{3}(\beta^{(2)}_{2})-20s.$$
The condition can be written as $R_{3}(\beta^{(2)}_{1})> 0$ and $R_{3}(\beta^{(2)}_{2})< 0$ or  $T_{1}(\beta^{(2)}_{1})<0$ and $T_{1}(\beta^{(2)}_{2})> 0$. Now, if the three roots of $R_{3}$ are interlaced with those of $P_{4}$, then $$20s\in{ ]\sup_{k\in\{1,3\}}R^{0}_{3}(\alpha^{(4)}_{k});\inf_{k\in\{2,4 \}}R^{0}_{3}(\alpha^{(4)}_{k})[}\subset]R^{0}_{3}(\beta^{(2)}_{2});R^{0}_{3}(\beta^{(2)}_{1})[$$
under the assumption $$\sup_{k\in\{1,3\}}R^{0}_{3}(\alpha^{(4)}_{k})<\inf_{k\in\{2,4 \}}R^{0}_{3}(\alpha^{(4)}_{k}).$$
\subsubsection{Discussion : our interval is reduced to a point}
Now for a better understanding of our assumption, we need to precise when
$$20s=\sup_{k\in\{1,3\}}R^{0}_{3}(\alpha^{(4)}_{k})=\inf_{k\in\{2,4 \}}R^{0}_{3}(\alpha^{(4)}_{k}).$$
This corresponds with two cases of multiple roots  for $P_{5}$: one special case of two real double roots and another one of a real triple root. This is the subject of the two following paragraphs.

\paragraph{Case of two real double roots}

The polynomial $P_5$ has two double roots $a$ and $b$ which are roots of $P_4$ and $R_3$. The Euclidean division of $ P_{5} $ by $ (x-a)^{2} (x-b)^2 $ gives a remainder that is identically zero if:
$$\left\{\begin{array}{l}
3(a^2+b^2)+4ab+\frac{10p}{3}=0\\
-2(a^3+b^3)-8ab(a+b)+10q=0\\
7a^2b^2+4ab(a^2+b^2)+20r=0\\
20s-2a^2b^2(a+b)=0.
\end{array}\right.
$$
We deduce
$$\left\{\begin{array}{l}
(ab)^2- \frac{8p}{3}ab+12r=0\\
(a+b)^3+\frac{2p}{3}(a+b)-2q=0\\
a+b=3q/[ab-\frac{2p}{3}]\\
-2(ab-\frac{2p}{3})^3+2p(ab-\frac{2p}{3})^2+27q^2=0.
\end{array}\right.
$$
\begin{enumerate}
\item The discriminant of the first equation in $(ab)$ is not negative if $r\leq \frac{4p^2}{27}$.
\item The discriminant of the second equation in $(a+b)$ must be not negative :  $-\frac{4}{27}(8p^3+729q^2)\geq 0$.
\item If $q\neq 0$, we conclude having two double real roots if $8p^3+729q^2\leq 0$ and $r\leq \frac{4p^2}{27}$.
\item The case $q=0$ gives  $r=\frac{5p^2}{36}$, $a=-b=\sqrt{\frac{-5p}{3}}$ and $s=0$ or $r=\frac{p^2}{9}$ and $ab=\frac{2p}{3}=-2\sqrt{r}$ and $s=\frac{2\sqrt{-2p}p^2}{45\sqrt{3}}$.
\item The third equation gives the particular case $ab=\frac{2p}{3}$: this implies that $q=0$, $r=p^2/9$ and $a+b=\pm \sqrt{-2p/3}$.
\end{enumerate}
\paragraph{Case of a real triple root} 
It also matches a double root of $P_{4}$ and $R_3$, a root of $P_{3}$ and $R_{2}$.
We can look for the order five under which conditions $ P_{5} $ has a triple root. If we divide $ P_{5} $ by $ (x-a)^3 $, we have, if the remainder is zero, that:
$$\left\{\begin{array}{l}
a^3+ap+q=0\\
4r-3a^4-2a^2p=0\\
6a^5+20s+\frac{10a^3p}{3}=0.
\end{array}\right.$$
The conditions on $a$ are those for a double root at the order four or quadruple at the order six (see later). The roots of  $4r-3a^4-2a^2p=0$ are
$$\pm\frac{1}{3}\sqrt{-3p\pm 3\sqrt{p^2+12r}}.$$
When $r=-\frac{p^2}{12}$ then $4p^3+27q^2=0$ and $ s=\frac{-pq}{30}$. One of these values (according with the sign of $q$) $$\sqrt{-p/3}, -\sqrt{-p/3}$$ should be a double root of $P_3$, triple root of $P_4$ and quadruple of $P_5$.\\
The quotient of $ P_{5} $ by $ (x-a)^3 $ gives
$$6a^2+3ax+x^2+\frac{10p}{3}.$$
So, if, in addition to the triple root, we also want two real roots: the discriminant of the previous polynomial is not negative if $$a\in\left]-\frac{2\sqrt{-2p}}{3};\frac{2\sqrt{-2p}}{3}\right[.$$
If $a=\pm 2\sqrt{-2p}/3$, $P_{5} $ has a triple root $ a $ and a double root.

\subsubsection{Discussion : our interval is empty or not empty}
Under the assumptions
$$p<0,\quad q\in \left]-2(-\frac{p}{3})^{3/2};2(-\frac{p}{3})^{3/2}\right[,\quad 4r\in ]R_{2}^{0}(\alpha_{2}^{(3)});\inf_{k\in{\{1,3\}}}R^{0}_{2}(\alpha_{k}^{(3)})[,$$
the assumption 
\begin{equation}\label{ass}
\sup_{k\in\{1,3\}}R^{0}_{3}(\alpha^{(4)}_{k})<\inf_{k\in\{2,4 \}}R^{0}_{3}(\alpha^{(4)}_{k})
\end{equation}
 may not be satisfied if
$R^{0}_{3}(\alpha^{(4)}_{1})>R^{0}_{3}(\alpha^{(4)}_{4})$ with $\alpha^{(4)}_{1}$ the smallest root and  $\alpha^{(4)}_{4}$ the biggest root. So, as the sum of the roots is zero, $\alpha^{(4)}_{1}<0$ and $\alpha^{(4)}_{4}>0$. Taking $X=\alpha^{(4)}_{1}\alpha^{(4)}_{4}$ (note, as the product of the four roots is equal to $4r$, $X\leq -2\sqrt{|r|}$) , the inequality (\ref{ass}) can be written as  
\begin{equation}\label{ass2}
(\frac{6r}{p}-p)X(X^2-4r)+(\frac{9q^2}{p}+2r)X^2-8r^2<0.
\end{equation}
\begin{enumerate}
\item $r<0$
\begin{enumerate}
\item This inequality is satisfied if $r<0$ or better if $r<\frac{-9q^2}{2p}$. $P_5$ has five distinct real roots.
\item If, for example, $\alpha^{(3)}_{1}$ or $\alpha^{(3)}_{3}$ according with the sign of $q$ is lower in absolute value than $|\frac{-3q}{p}|$ ($\inf_{k\in{\{1,3\}}}R^{0}_{2}(\alpha_{k}^{(3)})<0$), the proposed interval at order four for $r$ is included in $]-\infty; 0[$: $P_5$ has five distinct real roots.
\end{enumerate}
\item $r>0$
\begin{enumerate}
\item As $X\leq -2\sqrt{r}$, the inequality (\ref{ass2}) is satisfied for $r\leq \frac{p^2}{9}$: $P_5$ has five distinct real roots.
\item If $r>p^2/9$, we must look for negative roots of the polynomial (\ref{ass2}) which are less than $-2\sqrt{r}$.
\item If $r=\frac{4p^2}{27}$ and $729q^2+8p^3=0$, the polynomial with variable $X$ in (\ref{ass2}) has a not positive double root $\frac{4p}{3}$. By taking $X=4p/3$, $P_5$ has a real triple root and a real double (here, $s=4pq/5$).
\end{enumerate}
\end{enumerate}

  \subsection{Polynomial of  order six}
 In this section, another time we compare our assumptions with the Sturms' assumptions. Then, we explain how choosing $t$ and we describe some cases of multiple roots. Take

$$S_{6}(x)=P_{6}(x)=x^6+5px^4+20qx^3+60rx^2+120sx+120t$$
and
$$S_{5}(x)=P^{'}_{6}(x)=6P_{5}(x)=6x^5+20px^3+60qx^2+120rx+120s.$$
We call $\alpha^{(5)}_{1},\alpha^{(5)}_{2},\alpha^{(5)}_{3},\alpha^{(5)}_{4},\alpha^{(5)}_{5}$ the five distinct roots of $P_5$. We find:
$$S_{4}(x)=R_{4}(x)=-\left(\frac{5}{3}px^4+10qx^3+40rx^2+100sx+120t\right)$$
with $P_{6}(x)=\frac{x}{6}P^{'}_{6}(x)- S_{4}(x)$. We multiply $P^{'}_{6}$ by $\frac{p^2}{4}$ and divide by $S_{4}$. It provides $S_{3}$:
$$\begin{array}{lll}
 S_{3}(x)&=&-[(5p^3+54q^2-36pr)x^3+(216qr+15p^2q-90ps)]x^2\\
 &&-[(540qs+30p^2r-108pt)x+648qt+30p^2s].
 \end{array}$$
 Now, we multiply $S_{4}$ by $$\frac{(5p^3+54q^2-36pr)^{2}}{15p^2}.$$ We divide this polynomial by $S_{3}$ and we obtain as the opposite of the remainder:
$$S_{2}(x)=-a_{S_2}x^2-b_{S_2}x-c_{S_2}$$
with $$\left\{\begin{array}{lll}
a_{S_2}&=&-12(5p^3+54q^2-36pr)t-50p^4r-1620q^2pr+840p^2r^2+300p^2qs\\
&&-3456r^3+4320rqs+25p^3q^2+540q^4-900ps^2\\
b_{S_2}&=&-150p^4s-2520psq^2+2580p^2sr+180p^2qt-8640r^2s+2592rqt\\
&&+50qp^3r
+1080q^3r-1440qpr^2+5400s^2q-1080pst\\
c_{S_2}&=&-200tp^4-3240tpq^2+2880tp^2r-10368tr^2+50qp^3s+1080q^3s\\
&&-1440qprs+300p^2s^2+6480sqt
\end{array}\right.$$
Likewise:
$$S_{1}(x)=-a_{S_1}x-b_{S_1}$$
with $$\left\{\begin{array}{lll}
a_{S_1}&=&-135000s^4+250q^2p^3r^2+48600q^3spr+16200tp^2rq^2
-16200p^2s^2t\\&&+194400s^2rt-38880qt^2s-27900p^2s^2q^2+23100p^3s^2r
+60480tr^3p\\&&-250tp^4q^2+500tp^5r+900sqtp^3-45360sqtpr-3240p^2rt^2-69120r^5\\&&+3888pt^3-750p^4r^3+600p^4t^2-32400qp^2r^2s
-194400rq^2s^2-129600ps^2r^2\\&&-21600pq^2r^3
-38880tr^2q^2-500q^3sp^3-11100tr^2p^3-31104r^2t^2\\&&
+45360sq^3t-5400tpq^4+3240q^2t^2p+

241920r^3qs+144000ps^3q\\&&+1750qp^4rs+5400q^4r^2
-10800q^5s-1125p^5s^2+14400r^4p^2\\
b_{S_1}&=&-194400rq^2s^2-135000s^4+250q^2p^3r^2+48600q^3spr+16200tp^2rq^2\\&&
-16200p^2s^2t
+194400s^2rt-38880qt^2s-27900p^2s^2q^2+23100p^3s^2r\\&&+60480tr^3p-250tp^4q^2+500tp^5r+900sqtp^3-45360sqtpr-3240p^2rt^2\\&&-69120r^5
+3888pt^3-750p^4r^3+600p^4t^2-32400qp^2r^2s-
194400rq^2s^2\\&&-129600ps^2r^2-21600pq^2r^3-38880tr^2q^2
-500q^3sp^3-11100tr^2p^3-31104r^2t^2\\&&
+45360sq^3t-5400tpq^4+3240q^2t^2p+241920r^3qs+
144000ps^3q+1750qp^4rs\\&&+5400q^4r^2-10800q^5s-1125p^5s^2+14400r^4p^2
\end{array}\right.$$
Finally, $S_{0}$ is
$$ S_{0}=K_{1}\Delta(P_{6})=K_{2}\displaystyle\prod_{i=1}^{5}P_{6}(\alpha^{(5)}_{i})
=-K_{2}\displaystyle \prod_{i=1}^{5}R_{4}(\alpha^{(5)}_{i})=-K_{2}\displaystyle \prod_{i=1}^{5}(R^{0}_{4}(\alpha^{(5)}_{i})-120t)$$ 
According with Sturm's theorem, the assumptions for having six real distinct roots for the order six are:
$$\left\{\begin{array}{l}
p<0\\
5p^3+54q^2-36pr<0\\
a_{S_2}<0\\
a_{S_1}<0\\
S_0>0.
\end{array}\right.$$
In these cases, $a_{S_2}$ is a polynomial of order one in $t$ and $a_{S_1}$ a polynomial of order three in $t$. The study of the
sign will remove the unnecessary intervals of $\Delta(P_{6})>0$.\\
By our method, we must assume the conditions described in order five on $p,q,r,s$.

\paragraph{Choice of $t$:}
$R_{4}$ must have distinct roots. For this, we need three distinct roots $\beta^{(3)}_{1},\beta^{(3)}_{2},\beta^{(3)}_{3}$ of $R_{3}=R'_4/5$, $$R_{3}(x)=-\left(\frac{4}{3}px^3+6qx^2+16rx+20s\right).$$ We note once more
$$R_{4}^{0}(x)=-\left(\frac{5}{3}px^4+10qx^3+40rx^2+100sx\right).$$ It is necessary that
$R_{4}(\beta^{(3)}_{2})>0$ and that $R_{4}(\beta^{(3)}_{1})<0$ and $R_{4}(\beta^{(3)}_{3})<0$  
or else $R^{0}_{4}(\beta^{(3)}_{2})-120t>0$ and $R^{0}_{4}(\beta^{(3)}_{1})-120t<0$ and  $R^{0}_{4}(\beta^{(3)}_{3})-120t<0$.  Then we have to show that the four
roots of $R_{4}$ obtained are interlaced with $\alpha^{(5)}_{1},\alpha^{(5)}_{2},\alpha^{(5)}_{3},\alpha^{(5)}_{4},\alpha^{(5)}_{5}$ which will imply that
$120t$ belongs to the interval we want, that is to say,
$$120t\in{ ]\sup_{k\in\{2,4\}}R^{0}_{4}(\alpha^{(5)}_{k});\inf_{k\in\{1,3,5 \}}R^{0}_{4}(\alpha^{(5)}_{k})[}\subset ]\sup_{k\in\{1,3 \}}R^{0}_{4}(\beta^{(3)}_{k});R_{4}^{0}(\beta^{(3)}_{2})[ $$
under the assumption
$$\sup_{k\in\{2,4\}}R^{0}_{4}(\alpha^{(5)}_{k})<\inf_{k\in\{1,3,5 \}}R^{0}_{4}(\alpha^{(5)}_{k}).$$
To ensure that $R_{4}$ had four real roots, we divide $R_{4}$ by $R^{'}_{4}$. The remainder  will give a polynomial $T_{2}$ of degree two whose we change the sign
$$T_{2}(x)=\frac{5}{4}\frac{(-9q^2+16pr)x^2}{p}+\frac{15(5ps-2qr)x}{p}+120t-\frac{75qs}{2p}.$$ 
The fact that the roots of this polynomial $ \gamma ^{(2)} _ {1} $ and $\gamma^{(2)}_{2} $ which are a function of $ t $ are interlaced with
  the three distinct roots $ \beta ^ {(3)} _ {1}, \beta ^ {(3)} _ {2}, \beta ^ {(3)} _ {3} $ of $ R ^ { '} _ {4} $ give the condition for $ t $
  $$ \beta ^ {(3)} _ {1} <\gamma ^ {(2)} _ {1} <\beta ^ {(3)} _ {2} <\gamma ^ {(2)} _ {2} <\beta ^ {(3)} _ {3}. $$
More precisely, we consider $ T_ {2} $ whose coefficient of highest degree is positive: if $ p <0 $ and $ \Delta_ {2}> 0 $, so the conditions described in the order four give $ T_ {2} (\beta^{(3)}_{2}) <0 $, $ T_{2} (\beta ^ {(3)} _ {1})> 0 $ and $ T_{ 2} (\beta ^ {(3)} _ {3})> 0 $. This is equivalent to $ R_{4}(\beta ^ {(3)} _ {2})> 0 $, $ R_{4 } (\beta ^ {(3)}_{1}) <0 $ and $ R_{4} (\beta ^ {(3)}_{3}) <0 $; this interval for $ t $ is included in that which ensures that $ T_ {2} $ has two distinct real roots: $ T_{2}(\gamma ^{(1)}_{1}) <0 $ if $ \gamma ^ {(1)} _ {1} $ is the root of $ T_ {2} ^ {'} $ that we defined at the order five.
  \vspace{0.5cm}

Some bounds for these intervals can be obtain by studying different cases of multiple roots for the order six.
In that cases, several polynomials of the Sturm's sequence vanishes identically. Let us now describe three particular cases.
 
\paragraph{Multiple Roots}

We can focus for the order six under which conditions $ P_{6} $ has a root of multiplicity five. We already meet the case at the order three, four and five. If we divide $ P_{6} $ by $ (x-a)^5 $, we have, if the remainder is zero, that:
$$a^2=-\frac{p}{3}, \quad q^2=4a^6=-\frac{4p^3}{27}, \quad 4r=-3a^4=-\frac{p^2}{3},\,s=\pm \frac{p^2\sqrt{-p}}{45\sqrt{3}}, \quad t=-\frac{a^6}{24}=\frac{p^3}{648}.$$
In that case, we also have a quadruple root for $P_{5}$ and for $R_{4}$.
 Otherwise, if we want a real root of multiplicity four and a real double root, then $b=-2a$ with $a$ such that
 $$a^2=-\frac{5p}{6}, \quad q^2=\frac{a^6}{25}=-\frac{5p^3}{216}, \quad 4r=\frac{3a^4}{5}=\frac{5p^2}{12},\,
s=\pm \frac{p\sqrt{-5p}}{72\sqrt{6}}, \quad t=\frac{a^6}{30}=-\frac{25p^3}{1296}.$$
This case corresponds also with a triple root for $P_{5}$ and for $R_{4}$.\newline
At last, the case of two real triple roots for which  we find two double roots for $P_{5}$ and for $R_{4}$.\newline
 It comes:
 $$a^2=-\frac{5p}{3}, \quad q=0, \quad 4r=\frac{a^4}{5}=\frac{5p^2}{9},\,s=0, \quad t=-\frac{a^6}{120}=\frac{25p^3}{648}.$$
 In these three cases,  we remark that $S_{3}$ vanishes
identically.

\subsection{ Polynomial of order seven}
In this section, according with our theorem and our family of polynomials, we only explain the choice of the last constant.
In that case,we have
$$\begin{array}{l}
P_{7}(x)=x^7+7px^5+35qx^4+140rx^3+420sx^2+840tx+840u\\
P^{'}_{7}(x)=7(x^6+5px^4+20qx^3+60rx^2+120sx+120t).
\end{array}$$
We call $\alpha^{(6)}_{1},\alpha^{(6)}_{2},\alpha^{(6)}_{3},\alpha^{(6)}_{4},\alpha^{(6)}_{5},\alpha^{(6)}_{6}$ the six distinct real roots of $P^{'}_{7}$.
We define :
$$R_{5}(x)=-(2px^5+15qx^4+80rx^3+300sx^2+720tx+840u)$$
with
$$P_{7}(x)=\frac{x}{7}P^{'}_{7}(x)- R_{5}(x).$$
We put once again 
$$R^{0}_{5}(x)=-(2px^5+15qx^4+80rx^3+300sx^2+720tx).
$$
We assume the same assumptions for $p,q,r,s,t$ of the order six.

\paragraph{ Choice of $u$:}
As $R^{'}_{5}=6R_4$,  we call $\beta^{(4)}_{1},\beta^{(4)}_{2},\beta^{(4)}_{3},\beta^{(4)}_{4}$ their four distinct real roots. Then, the Euclidean division of $R_{5}$ by  $R^{'}_{5}$ gives changing the sign of the remainder:
$$T_{3}(x)=\frac{2}{p}(-9q^2+16pr)x^3+\frac{36(5ps-2qr)}{p}x^2$$
$$+\frac{36s(16pt-5qs)}{p}x+840u-\frac{216qt}{p}.$$ 
The Euclidean division of $T_{3}$ by $T^{'}_{3}$ gives changing the sign of the remainder a polynomial $U_{1}$ of degree one whose numerator of the term in $x $ is none other than the discriminant of the polynomial $ T_{2} $ of the order six which is positive to have two real roots $(T_{2}(\gamma^{(1)}_{1})<0)$. 
We deduce the following expression for $U_{1}$:
$$U_{1}(x)=-\left[\frac{24(-144ptq^2+256p^2tr+45q^3s+40qspr-150p^2s^2-24q^2r^2)}{p(-9q^2+16pr)}x\right.$$
$$\left.+\frac{24(81q^3t-48qtpr-315puq^2+560p^2ur-240p^2st+75ps^2q-30q^2rs)}{p(-9q^2+16pr)}\right]$$
whose the root $\delta^{(1)}_{1}$ is a polynomial of degree one in $u$:
 $$\delta^{(1)}_{1}=\frac{-[81q^3t-48qtpr-35pu(9q^2-16pr)-240p^2st+75ps^2q-30q^2rs]}{-16pt(9q^2-16pr)+45q^3s+40qspr-150p^2s^2-24q^2r^2}.$$
If this root is interlaced with those $ \gamma^{(2)} _ {1} $ and $\gamma^{(2)}_{2}$ of $ T_{2} $ - as established for the  degree three-, $ T_{3} $ has three real distinct roots. If these roots are interlaced with $\beta^{(4)}_{1},\beta^{(4)}_{2},\beta^{(4)}_{3},\beta^{(4)}_{4}$, this gives an interval for $ u $ included in the previous one. So according with the result for the order five, $ R_{5} $ has five real distinct roots.
Finally, if these five real roots are interlaced with
$\alpha^{(6)}_{1},\alpha^{(6)}_{2},\alpha^{(6)}_{3},\alpha^{(6)}_{4},\alpha^{(6)}_{5},\alpha^{(6)}_{6}$, we conclude that $P_{7}$ has seven real distinct roots. We find:
 $$840u\in{ \left]\sup_{k\in\{1,3,5\}}R^{0}_{5}(\alpha^{(6)}_{k}),\inf_{k\in\{2,4,6 \}}R^{0}_{5}(\alpha^{(6)}_{k})\right[},$$
if we have $$\sup_{k\in\{1,3,5\}}R^{0}_{5}(\alpha^{(6)}_{k})<\inf_{k\in\{2,4,6 \}}R^{0}_{5}(\alpha^{(6)}_{k}).$$
\vspace{0.5cm}

\section{Concluding Remarks}

Our assumptions are explicit and depend on the roots of the previous order. That is why after the order five, things become harder. Indeed, as we recall in the introduction, exact expressions of the roots are unknown. Obviously, for cases of order three and four, Cardano's, Ferrari's,  Descartes's or Euler's formula of the roots are available and can be used for the order five. Notice that, for example  in  different cases of multiple roots, the resolvant cubic takes a nice and simple form.\\ 
Of course, recall that there exists methods, for example, Brings-Jerrard efficient for solving polynomials of degree five. Different transformations are needed in order to obtain expressions of the roots.\newline
Otherwise,  Cayley  \cite{Cayley61} and more  recently  \cite{Dum91} and \cite{Laval05} proposed different methods for some kind of polynomials of degree six those roots as functions of the roots for solvable quintics. On the other side, Ramanujan \cite{Bern94} solve some polynomials of degree three, four, five, six and seven. His approach is very original. It seems that he often started with roots having product one. In this paper, we rather assume that the sum of the roots is zero even if our result is yet true in the general case, but it takes a more complex form. Naturally, this is particularly true for the Sturm's sequence. Some connections with the theory of elliptic functions would be surely promising and successful.

\end{document}